\documentclass[reqno,a4paper]{amsart}
\usepackage{amssymb}%
\usepackage[left=2cm,right=2cm]{geometry}%
\usepackage{ifpdf}
\usepackage{listings}
\usepackage{color,xcolor}
\usepackage{colortbl}
\usepackage{multirow}
\usepackage{booktabs}
\usepackage{float}
\usepackage{graphicx}
\usepackage{algorithmic}
\usepackage[linesnumbered,ruled,vlined]{algorithm2e}
\usepackage{subfig}
\usepackage{epstopdf}
\usepackage[colorlinks,linkcolor=red,anchorcolor=red,citecolor=red]{hyperref}
\usepackage[numbers,sort&compress]{natbib}
\usepackage{enumitem}

\theoremstyle{plain}
\newtheorem{thm}{Theorem}[section]

\theoremstyle{remark}
\newtheorem{remark}{Remark}[section]
\theoremstyle{definition}

\numberwithin{equation}{section}
\allowdisplaybreaks[4]

\begin{document}
\title[]{Recurrence Relations for the Maclaurin Coefficients of Products of Elementary Functions and the Bessel Function}
\author{Zhong-Xuan Mao, Jing-Feng Tian*}

\address{Zhong-Xuan Mao \\
Hebei Key Laboratory of Physics and Energy Technology\\
Department of Mathematics and Physics\\
North China Electric Power University \\
Yonghua Street 619, 071003, Baoding, P. R. China}
\email{maozhongxuan000\symbol{64}gmail.com}

\address{Jing-Feng Tian\\
Hebei Key Laboratory of Physics and Energy Technology\\
Department of Mathematics and Physics\\
North China Electric Power University\\
Yonghua Street 619, 071003, Baoding, P. R. China}
\email{tianjf\symbol{64}ncepu.edu.cn}

\begin{abstract}
In this paper, we investigate recurrence relations for the Maclaurin coefficients of the Products of a elementary function and the Bessel function of the first kind $\mathcal{J}(z) = h(z) J_\nu(z)$ and the modified Bessel function of the first kind
$\mathcal{I}(z) = h(z) I_\nu(z)$ in the complex plane corresponding to several specific choices of $h(z)$. In particular, we specialize $h(z)$ as
$e^{pz}$, $(1-\theta z)^p$, $e^{-p \arctan z}$, $\sin(pz)$, $\cos(pz)$, $\sinh(pz)$, $\cosh(pz)$, $\arcsin(pz)$ and $\arccos(pz)$.

\end{abstract}

\thanks{*Corresponding author: Jing-Feng Tian(tianjf\symbol{64}ncepu.edu.cn)}

\maketitle


\section{Introduction}

The Bessel function of the first kind is defined by the series expansion
\begin{equation} \label{J-series}
J_\nu(z) = \sum_{n=0}^\infty \frac{(-1)^n}{n! \Gamma(\nu+n+1)} \Big(\frac{z}{2}\Big)^{2n+\nu}, \quad z\in\mathbb{C}, \quad \nu \notin \{-1,-2,\cdots\},
\end{equation}
and the modified Bessel function of the first kind is defined by
\begin{equation} \label{I-series}
I_\nu(z) = \sum_{n=0}^\infty \frac{1}{n! \Gamma(\nu+n+1)} \Big(\frac{z}{2}\Big)^{2n+\nu}, \quad z\in\mathbb{C}, \quad \nu \notin \{-1,-2,\cdots\}.
\end{equation}
These functions are not only of fundamental theoretical interest \cite{Bowman-1958, Luke-1962} but also find extensive applications across a wide range of scientific and engineering disciplines \cite{Korenev-2002}.

In the theory of special functions, one frequently encounters expressions consisting of the product of an elementary function---such as an exponential, logarithmic, trigonometric, or hyperbolic function---and a Bessel function (or a modified Bessel function). The analysis of such products often presents computational and symbolic challenges that are not present when dealing with isolated Bessel functions. Motivated by these observations, we systematically investigate the series expansions of this class of composite functions in the present paper.

Specifically, let $\mathcal{J}(z)=h(z)J_\nu(z)$ and $\mathcal{I}(z)=h(z)I_\nu(z)$. We aim to establish explicit recurrence relations for the coefficients of their Maclaurin series in the complex plane. We consider a broad spectrum of specific choices for the auxiliary function $h(z)$. In particular, we specialize $h(z)$ to include:
\begin{itemize}
    \item Exponential and algebraic forms: $e^{pz}$, $(1-\theta z)^p$, and the composite form $e^{-p\arctan z}$;
    \item Trigonometric and hyperbolic functions: $\sin(pz)$, $\cos(pz)$, $\sinh(pz)$, and $\cosh(pz)$;
    \item Inverse trigonometric functions: $\arcsin(pz)$ and $\arccos(pz)$.
\end{itemize}

It is worth noting that the method of deriving recurrence relations for expansion coefficients has proven to be a powerful tool in the study of other special functions. Similar techniques have been successfully applied to complete elliptic integrals \cite{Yang-JMAA-2018} and Gaussian hypergeometric functions \cite{Chen-RACSAM-2021, Yang-PAMS-2025}. These results have served as effective analytical instruments in subsequent research, such as establishing the absolute monotonicity of specific functions \cite{Wu-IJPAM-2025}. Moreover, Yang's 2025 internal report \cite{Yang2025Report} offered an early systematic treatment of recurrence methods for special-function expansions, laying a foundation for subsequent studies.

\section{Main results}

\subsection{Recurrence relations for the Maclaurin coefficients of products of elementary functions and the Bessel function of the first kind}

In this part, we investigate the recurrence relations for the Maclaurin coefficients of functions of the form $\mathcal{J}(z) = h(z) J_\nu(z)$, where $h(z)$ is chosen as certain specific elementary functions.

We first consider the product of the exponential function and the Bessel function of the first kind.

\begin{thm} \label{thm-exp-J}
Let $\nu,p \in \mathbb{C}$ and $\nu \notin \{-\frac{1}{2},-\frac{3}{2},\cdots\}$. Suppose that
\begin{equation*}
e^{pz} J_\nu(z) = \frac{1}{2^\nu \Gamma (\nu +1)} \sum_{n=0}^\infty u_n z^{n+\nu}, \quad z \in \mathbb{C}.
\end{equation*}
Then $u_0 = 1$, $u_1 = p$ and
\begin{equation*}
u_{n+1} = \frac{p (2 \nu +2 n+1)}{(n+1) (2 \nu +n+1)} u_n -\frac{p^2+1}{(n+1) (2 \nu +n+1)} u_{n-1}.
\end{equation*}
\end{thm}

Noting that $\sinh(pz) = (e^{pz} - e^{-pz})/2$ and $\cosh(pz) = (e^{pz} + e^{-pz})/2$, we have
\begin{equation*}
\sinh(pz) J_\nu(z) = \frac{1}{2} \Big( e^{pz} J_\nu(z) - e^{-pz} J_\nu(z) \Big), \end{equation*}
and
\begin{equation*}
\cosh(pz) J_\nu(z) = \frac{1}{2} \Big( e^{pz} J_\nu(z) + e^{-pz} J_\nu(z) \Big),
\end{equation*}
which lead to the following two theorems from Theorem \ref{thm-exp-J}.

\begin{thm} \label{thm-sinh-J-2}
Let $\nu,p \in \mathbb{C}$ and $\nu \notin \{-\frac{1}{2},-\frac{3}{2},\cdots\}$. Suppose that
\begin{equation*}
\sinh(pz) J_\nu(z) = \frac{1}{2^{\nu+1} \Gamma (\nu +1)} \sum_{n=0}^\infty (u_n - v_n) z^{n+\nu}, \quad z \in \mathbb{C}.
\end{equation*}
Then $u_0 = 1$, $u_1 = p$, $v_0 = 1$, $v_1 = -p$,
\begin{equation*}
u_{n+1} = \frac{p (2 \nu +2 n+1)}{(n+1) (2 \nu +n+1)} u_n -\frac{p^2+1}{(n+1) (2 \nu +n+1)} u_{n-1},
\end{equation*}
and
\begin{equation*}
v_{n+1} = - \frac{p (2 \nu +2 n+1)}{(n+1) (2 \nu +n+1)} v_n -\frac{p^2+1}{(n+1) (2 \nu +n+1)} v_{n-1}.
\end{equation*}
\end{thm}

\begin{thm} \label{thm-cosh-J-2}
Let $\nu,p \in \mathbb{C}$ and $\nu \notin \{-\frac{1}{2},-\frac{3}{2},\cdots\}$. Suppose that
\begin{equation*}
\cosh(pz) J_\nu(z) = \frac{1}{2^{\nu+1} \Gamma (\nu +1)} \sum_{n=0}^\infty (u_n + v_n) z^{n+\nu}, \quad z \in \mathbb{C}.
\end{equation*}
Then $u_0 = 1$, $u_1 = p$, $v_0 = 1$, $v_1 = -p$,
\begin{equation*}
u_{n+1} = \frac{p (2 \nu +2 n+1)}{(n+1) (2 \nu +n+1)} u_n -\frac{p^2+1}{(n+1) (2 \nu +n+1)} u_{n-1},
\end{equation*}
and
\begin{equation*}
v_{n+1} = - \frac{p (2 \nu +2 n+1)}{(n+1) (2 \nu +n+1)} v_n -\frac{p^2+1}{(n+1) (2 \nu +n+1)} v_{n-1}.
\end{equation*}
\end{thm}

Noting that $\sin(pz) = \frac{e^{ipz} - e^{-ipz}}{2i}$ and $\cos(pz) = \frac{e^{ipz} + e^{-ipz}}{2}$, we have
\begin{equation*}
\sin(pz) J_\nu(z) = \frac{1}{2i} \Big( e^{ipz} J_\nu(z) - e^{-ipz} J_\nu(z) \Big), \end{equation*}
and
\begin{equation*}
\cos(pz) J_\nu(z) = \frac{1}{2} \Big( e^{ipz} J_\nu(z) + e^{-ipz} J_\nu(z) \Big),
\end{equation*}
which lead to the following two theorems from Theorem \ref{thm-exp-J}.

\begin{thm} \label{thm-sinh-J-2}
Let $\nu,p \in \mathbb{C}$ and $\nu \notin \{-\frac{1}{2},-\frac{3}{2},\cdots\}$. Suppose that
\begin{equation*}
\sinh(pz) J_\nu(z) = \frac{1}{2^{\nu+1} i \Gamma (\nu +1)} \sum_{n=0}^\infty (u_n - v_n) z^{n+\nu}, \quad z \in \mathbb{C}.
\end{equation*}
Then $u_0 = 1$, $u_1 = ip$, $v_0 = 1$, $v_1 = -ip$,
\begin{equation*}
u_{n+1} = \frac{ip (2 \nu +2 n+1)}{(n+1) (2 \nu +n+1)} u_n + \frac{p^2-1}{(n+1) (2 \nu +n+1)} u_{n-1},
\end{equation*}
and
\begin{equation*}
v_{n+1} = - \frac{ip (2 \nu +2 n+1)}{(n+1) (2 \nu +n+1)} v_n + \frac{p^2-1}{(n+1) (2 \nu +n+1)} v_{n-1}.
\end{equation*}
\end{thm}

\begin{thm} \label{thm-cosh-J-2}
Let $\nu,p \in \mathbb{C}$ and $\nu \notin \{-\frac{1}{2},-\frac{3}{2},\cdots\}$. Suppose that
\begin{equation*}
\cosh(pz) J_\nu(z) = \frac{1}{2^{\nu+1} \Gamma (\nu +1)} \sum_{n=0}^\infty (u_n + v_n) z^{n+\nu}, \quad z \in \mathbb{C}.
\end{equation*}
Then $u_0 = 1$, $u_1 = ip$, $v_0 = 1$, $v_1 = -ip$,
\begin{equation*}
u_{n+1} = \frac{ip (2 \nu +2 n+1)}{(n+1) (2 \nu +n+1)} u_n + \frac{p^2-1}{(n+1) (2 \nu +n+1)} u_{n-1},
\end{equation*}
and
\begin{equation*}
v_{n+1} = - \frac{ip (2 \nu +2 n+1)}{(n+1) (2 \nu +n+1)} v_n + \frac{p^2-1}{(n+1) (2 \nu +n+1)} v_{n-1}.
\end{equation*}
\end{thm}

Next, we consider the product of the Bessel function of the first kind with the factor $(1-\theta z)^p$.

\begin{thm}
Let $\nu,p,\theta,z \in \mathbb{C}$ and $\nu \notin \{-\frac{1}{2},-\frac{3}{2},\cdots\}$. Suppose that
\begin{equation*}
(1-\theta z)^p J_\nu(z) = \frac{1}{2^\nu \Gamma (\nu +1)} \sum_{n=0}^\infty u_n z^{n+\nu}, \quad |z| < \frac{1}{\theta}.
\end{equation*}
Then $u_0 = 1$, $u_1 = -\theta p$, $u_2 = \frac{1}{2} \theta ^2 (p-1) p-\frac{1}{4 (\nu +1)}$,
$u_3 = \frac{\theta  p}{4 (\nu +1)}-\frac{1}{6} \theta ^3 (p-2) (p-1) p$ and
\begin{equation*}
u_{n+1} = \beta_0(n) u_n + \beta_1(n) u_{n-1} + \beta_2(n) u_{n-2} + \beta_3(n) u_{n-3},
\end{equation*}
where
\begin{equation*}
\begin{aligned}
\beta_0(n) &= \frac{\theta  \left(2 n^2+4 \nu  n-2 n p-2 \nu  p-p\right)}{(n+1) (2 \nu +n+1)},\\
\beta_1(n) &= \frac{\theta ^2 (n-p-1) (2 \nu +n-p-1)+1}{(n+1) (2 \nu +n+1)},\\
\beta_2(n) &= \frac{2 \theta }{(n+1) (2 \nu +n+1)}, \\
\beta_3(n) &= -\frac{\theta ^2}{(n+1) (2 \nu +n+1)}.
\end{aligned}
\end{equation*}
\end{thm}

Next, we consider the product of $\exp(-p \arctan z)$ and the Bessel function of the first kind.

\begin{thm} \label{thm-exp-arctan-J}
Let $\nu,p \in \mathbb{C}$ and $\nu \notin \{-\frac{1}{2},-\frac{3}{2},\cdots\}$. Suppose that
\begin{equation*}
e^{-p \arctan z} J_\nu(z) = \frac{1}{2^\nu \Gamma (\nu +1)} \sum_{n=0}^\infty u_n z^{n+\nu}, \quad z \in \mathbb{C}.
\end{equation*}
Then $u_0 = 1$, $u_1 = - p$, $u_2 = \frac{p^2}{2}-\frac{1}{4 (\nu +1)}$,
$u_3 = \frac{1}{3} (p-\frac{p^3}{2})+\frac{p}{4 (\nu +1)}$,
$u_4 = \frac{1}{96} (\frac{3}{\nu ^2+3 \nu +2}+4 p^2 (-\frac{3}{\nu +1}+p^2-8))$ and
\begin{equation*}
u_{n+1} = \beta_0(n) u_n + \beta_1(n) u_{n-1} + \beta_2(n) u_{n-2} + \beta_3(n) u_{n-3} + \beta_4(n) u_{n-4},
\end{equation*}
where
\begin{equation*}
\begin{aligned}
\beta_0(n) &= -\frac{p (2 \nu +2 n+1)}{(n+1) (2 \nu +n+1)},\\
\beta_1(n) &= -\frac{-4 \nu +2 n (2 \nu +n-2)+p^2+3}{(n+1) (2 \nu +n+1)},\\
\beta_2(n) &= \frac{p (-2 \nu -2 n+5)}{(n+1) (2 \nu +n+1)}, \\
\beta_3(n) &= \frac{6 \nu -n (2 \nu +n-6)-11}{(n+1) (2 \nu +n+1)}, \\
\beta_4(n) &= -\frac{1}{(n+1) (2 \nu +n+1)}.
\end{aligned}
\end{equation*}
\end{thm}

Next, the following two theorems consider the products of trigonometric functions and the Bessel function of the first kind, and characterize their coefficients by a single recurrence relation.

\begin{thm} \label{thm-sin-J}
Let $\nu,p \in \mathbb{C}$ and $\nu \notin \{-\frac{1}{2},-\frac{3}{2},\cdots\}$. Suppose that
\begin{equation*}
\sin(pz) J_\nu(z) = \frac{1}{2^\nu \Gamma (\nu +1)} \sum_{n=0}^\infty u_n z^{n+\nu}, \quad z \in \mathbb{C}.
\end{equation*}
Then $u_{2k}=0$ for all $k\geq 0$, $u_1 = p$, $u_3 = -\frac{p^3}{6}-\frac{p}{4 (\nu +1)}$, $u_5 = \frac{p^5}{120}+\frac{p^3}{24 (\nu +1)}+\frac{p}{32 (\nu +1) (\nu +2)}$ and
\begin{equation*}
u_{2n+2} = \beta_0(n) u_{2n} + \beta_1(n) u_{2n-2} + \beta_2(n) u_{2n-4},
\end{equation*}
where
\begin{equation*}
\begin{aligned}
\beta_0(n) &= \frac{ \left(
\begin{aligned}
& -4 (p^2-1)n^4 - 16 (\nu -1)(p^2-1) n^3 \\
& \quad +2 \left(4 (\nu -6) \nu +(-12 (\nu -2) \nu -7) p^2+9\right) n^2 \\
& \quad +2 \left(-2 \nu  (2 \nu  (2 \nu +5)-9)+\left(3-2 \nu  \left(4 \nu ^2-2 \nu +7\right)\right) p^2+1\right) n \\
& \quad -(2 \nu -1) (2 \nu +1)^2 \left(2 (\nu +1) p^2+1\right)
\end{aligned}
\right)
}{\left(4 \nu ^2-1\right) (n+1) (n+2) (2 \nu +n+1) (2 \nu +n+2)},\\
\beta_1(n) &= -\frac{8 n^2 \left(p^4-1\right)+16 (\nu -2) n \left(p^4-1\right)+\left(p^2-1\right) \left(-4 \nu  (\nu +8)+(4 \nu  (5 \nu -8)+27) p^2+33\right)}{\left(4 \nu ^2-1\right) (n+1) (n+2) (2 \nu +n+1) (2 \nu +n+2)},\\
\beta_2(n) &= -\frac{4 \left(p^2-1\right)^3}{\left(4 \nu ^2-1\right) (n+1) (n+2) (2 \nu +n+1) (2 \nu +n+2)}.
\end{aligned}
\end{equation*}
\end{thm}

\begin{thm} \label{thm-cos-J}
Let $\nu,p \in \mathbb{C}$ and $\nu \notin \{-\frac{1}{2},-\frac{3}{2},\cdots\}$. Suppose that
\begin{equation*}
\cos(pz) J_\nu(z) = \frac{1}{2^\nu \Gamma (\nu +1)} \sum_{n=0}^\infty u_n z^{n+\nu}, \quad z \in \mathbb{C}.
\end{equation*}
Then $u_{2k-1}=0$ for all $k\geq 0$, $u_0 = 1$, $u_2 =-\frac{1}{4 (\nu +1)}-\frac{p^2}{2}$, $u_4 = \frac{1}{32 (\nu +1) (\nu +2)}+\frac{p^4}{24}+\frac{p^2}{8 (\nu +1)}$ and
\begin{equation*}
u_{2n+2} = \beta_0(n) u_{2n} + \beta_1(n) u_{2n-2} + \beta_2(n) u_{2n-4},
\end{equation*}
where $\beta_i(n)$ ($n=0,1,2$) is defined in Theorem \ref{thm-sin-J}.
\end{thm}

Next, we consider the product of hyperbolic functions and the Bessel function of the first kind.
\begin{thm} \label{thm-sinh-J}
Let $\nu,p \in \mathbb{C}$ and $\nu \notin \{-\frac{1}{2},-\frac{3}{2},\cdots\}$. Suppose that
\begin{equation*}
\sinh(pz) J_\nu(z) = \frac{1}{2^\nu \Gamma (\nu +1)} \sum_{n=0}^\infty u_n z^{n+\nu}, \quad z \in \mathbb{C}.
\end{equation*}
Then $u_{2k}=0$ for all $k\geq 0$, $u_1 = p$, $u_3 = \frac{p^3}{6}-\frac{p}{4 (\nu +1)}$, $u_5 = \frac{p^5}{120}-\frac{p^3}{24 (\nu +1)}+\frac{p}{32 (\nu +1) (\nu +2)}$ and
\begin{equation*}
u_{2n+2} = \beta_0(n) u_{2n} + \beta_1(n) u_{2n-2} + \beta_2(n) u_{2n-4},
\end{equation*}
where
\begin{equation*}
\begin{aligned}
\beta_0(n) &= \frac{
\left(
\begin{aligned}
& 4 \left(p^2+1\right)n^4 +16 (\nu -1)\left(p^2+1\right) n^3 \\
& \quad +2 \left(4 (\nu -6) \nu +(12 (\nu -2) \nu +7) p^2+9\right)n^2 \\
& \quad +2\left(-2 \nu  (2 \nu  (2 \nu +5)-9)+\left(2 \nu  \left(4 \nu ^2-2 \nu +7\right)-3\right) p^2+1\right) n \\
& \quad +(2 \nu -1) (2 \nu +1)^2 \left(2 (\nu +1) p^2-1\right)
\end{aligned}\right)
}{\left(4 \nu ^2-1\right) (n+1) (n+2) (2 \nu +n+1) (2 \nu +n+2)},\\
\beta_1(n) &= -\frac{8 \left(p^4-1\right) n^2+16 \left(p^4-1\right) (\nu -2) n+\left(p^2+1\right) \left((4 \nu  (5 \nu -8)+27) p^2+4 \nu  (\nu +8)-33\right)}{(n+1) (n+2) (n+2 \nu +1) (n+2 \nu +2) \left(4 \nu ^2-1\right)},\\
\beta_2(n) &= \frac{4 \left(p^2+1\right)^3}{\left(4 \nu ^2-1\right) (n+1) (n+2) (2 \nu +n+1) (2 \nu +n+2)}.
\end{aligned}
\end{equation*}
\end{thm}

\begin{thm}
Let $\nu,p \in \mathbb{C}$ and $\nu \notin \{-\frac{1}{2},-\frac{3}{2},\cdots\}$. Suppose that
\begin{equation*}
\cosh(pz) J_\nu(z) = \frac{1}{2^\nu \Gamma (\nu +1)} \sum_{n=0}^\infty u_n z^{n+\nu}, \quad z \in \mathbb{C}.
\end{equation*}
Then $u_{2k-1}=0$ for all $k\geq 0$, $u_0 = 1$, $u_2 = \frac{p^2}{2}-\frac{1}{4 (\nu +1)}$, $u_4 = \frac{1}{32 (\nu +1) (\nu +2)}+\frac{p^4}{24}-\frac{p^2}{8 (\nu +1)}$ and
\begin{equation*}
u_{2n+2} = \beta_0(n) u_{2n} + \beta_1(n) u_{2n-2} + \beta_2(n) u_{2n-4},
\end{equation*}
where $\beta_i(n)$ ($n=0,1,2$) is defined in Theorem \ref{thm-sinh-J}.
\end{thm}

At the end, we consider the product of an inverse trigonometric function and a Bessel function of the first kind.

\begin{thm} \label{thm-arcsin-J}
Let $\nu,p \in \mathbb{C}$ and $\nu \notin \{-\frac{1}{2},-\frac{3}{2},\cdots\}$. Suppose that
\begin{equation*}
\arcsin(pz) J_\nu(z) = \frac{1}{2^\nu \Gamma (\nu +1)} \sum_{n=0}^\infty u_n z^{n+\nu}, \quad z \in \mathbb{C}.
\end{equation*}
Then
$u_0 = u_2 = u_4 = u_6 = u_8 = u_{10} = u_{12} = 0$,
\begin{equation*}
u_1 = p, u_3 = \frac{p^3}{6}-\frac{p}{4 (\nu +1)}, u_5 = \frac{3 p^5}{40}-\frac{p^3}{24 (\nu +1)}+\frac{p}{32 (\nu +1) (\nu +2)},
\end{equation*}
\begin{equation*}
u_7 = \frac{5 p^7}{112}-\frac{3 p^5}{160 (\nu +1)}+\frac{p^3}{192 (\nu +1) (\nu +2)}-\frac{p}{384 (\nu +1) (\nu +2) (\nu +3)},
\end{equation*}
\begin{equation*}
u_9 = \frac{35 p^9}{1152}-\frac{5 p^7}{448 (\nu +1)}+\frac{3 p^5}{1280 (\nu +1) (\nu +2)}-\frac{p^3}{2304 (\nu +1) (\nu +2) (\nu +3)}+\frac{p}{6144 (\nu +1) (\nu +2) (\nu +3) (\nu +4)},
\end{equation*}
\begin{equation*}
\begin{aligned}
u_{11} & = \frac{63 p^{11}}{2816}-\frac{35 p^9}{4608 (\nu +1)}+\frac{5 p^7}{3584 (\nu +1) (\nu +2)}-\frac{p^5}{5120 (\nu +1) (\nu +2) (\nu +3)} \\
& \quad +\frac{p^3}{36864 (\nu +1) (\nu +2) (\nu +3) (\nu +4)}-\frac{p}{122880 (\nu +1) (\nu +2) (\nu +3) (\nu +4) (\nu +5)},
\end{aligned}
\end{equation*}
\begin{equation*}
\begin{aligned}
u_{13} & = \frac{231 p^{13}}{13312}-\frac{63 p^{11}}{11264 (\nu +1)}+\frac{35 p^9}{36864 (\nu +1) (\nu +2)}-\frac{5 p^7}{43008 (\nu +1) (\nu +2) (\nu +3)}\\
& \quad +\frac{p^5}{81920 (\nu +1) (\nu +2) (\nu +3) (\nu +4)}-\frac{p^3}{737280 (\nu +1) (\nu +2) (\nu +3) (\nu +4) (\nu +5)}\\
& \quad+\frac{p}{2949120 (\nu +1) (\nu +2) (\nu +3) (\nu +4) (\nu +5) (\nu +6)},
\end{aligned}
\end{equation*}
and
\begin{equation*}
u_{n+1} = \sum_{i=0}^{13} \beta_i(n) u_{n-i}, n\geq 13,
\end{equation*}
where
\begin{equation*}
\beta_0(n) = \frac{2 p^2 (\nu +n)}{(n+1) (2 \nu +n+1)},
\end{equation*}
\begin{equation*}
\beta_1(n) = - \frac{\left(
\begin{aligned}
& -(n-4) (n-3) (n-2) (n-1) \left(p^4+4 \left(4 \nu ^2-1\right) p^2+4\right) \\
& \quad -2 (2 \nu +1) (n-3) (n-2) (n-1) \left(p^4+8 \left(2 \nu ^2+\nu -1\right) p^2+4\right) \\
& \quad -2 (n-2) (n-1)  \left(4 \nu ^2-1\right) \left(4 (\nu +1) (2 \nu +1) p^2+1\right) \\
& \quad + (2 \nu -1) (2 \nu +1)^2+2 (2 \nu -1) (2 \nu +1) (n-1) \left(\nu  \left(p^4+2\right)+p^4+4\right)
\end{aligned}
\right)
}{\left(4 \nu ^2-1\right) n (n+1) (2 \nu +n) (2 \nu +n+1)},
\end{equation*}
\begin{equation*}
\beta_2(n) = - \frac{p^2 \left(
\begin{aligned}
& 4 (n-5) (n-4) (n-3) (n-2) p^2 + 2 (n-4) (n-3) (n-2) \left(p^4+(4 \nu  (3 \nu +2)-1) p^2+4\right) \\
& \quad +3 (n-3) (n-2)(2 \nu +1) \left(p^2 \left(6 \nu  (2 \nu +1)+p^2-8\right)+4\right) \\
& \quad +\left(4 \nu ^2-1\right) (n-2) \left(p^4+2 (\nu  (6 \nu +5)-2) p^2+2\right) \\
& \quad -2 (\nu +2) (2 \nu -1) (2 \nu +1)
\end{aligned}
\right)
}{\left(4 \nu ^2-1\right) n (n+1) (2 \nu +n) (2 \nu +n+1)},
\end{equation*}
\begin{equation*}
\beta_3(n) = - \frac{\left(
\begin{aligned}
& 2 (n-6) (n-5) (n-4) (n-3) p^2 \left(p^4+3 \left(4 \nu ^2-1\right) p^2+8\right)\\
& \quad 8 (n-5) (n-4) (n-3) p^2 \left(8 \nu +p^2 \left(\nu  \left(12 \nu  (\nu +1)+p^2-3\right)-3\right)+4\right) \\
& \quad -(n-4) (n-3) \left(p^8+12 \nu  p^6-6 \left(2 \nu  \left(2 \nu  \left(4 \nu ^2+6 \nu +1\right)-3\right)-3\right) p^4+\left(8-32 \nu ^2\right) p^2+8\right) \\
& \quad -(2 \nu +1) (n-3) \left(p^8+4 (2 \nu  (\nu +1)-3) p^6+6 p^4+16 (\nu +2) (2 \nu -1) p^2+8\right) \\
& \quad -(2 \nu -1) (2 \nu +1) \left(-p^4+(8 \nu +4) p^2-1\right)
\end{aligned}
\right)
}{\left(4 \nu ^2-1\right) n (n+1) (2 \nu +n) (2 \nu +n+1)},
\end{equation*}
\begin{equation*}
\beta_4(n) = - \frac{p^2 \left(
\begin{aligned}
& -8 (n-7) (n-6) (n-5) (n-4) p^4 \\
& \quad -2 (n-6) (n-5) (n-4) p^2 \left(p^4+(4 \nu  (3 \nu +4)+9) p^2+12\right) \\
& \quad -(n-5) (n-4) p^2 \left(72 \nu +(6 \nu +1) p^4+6 \left(12 (\nu +1) \nu ^2+\nu -1\right) p^2+12\right) \\
& \quad -(n-4) \left((2 \nu -3) (2 \nu +1) p^6+(2 \nu  (2 \nu  (2 \nu  (6 \nu +1)-9)-1)+4) p^4+(24 (\nu -2) \nu -26) p^2-8\right) \\
& \quad +(2 \nu +1) \left(p^4+2 (\nu  (6 \nu +5)-5) p^2+4\right)
\end{aligned}
\right)
}{\left(4 \nu ^2-1\right) n (n+1) (2 \nu +n) (2 \nu +n+1)},
\end{equation*}
\begin{equation*}
\beta_5(n) = - \frac{\left(
\begin{aligned}
& -(n-8) (n-7) (n-6) (n-5) p^4 \left(p^4+4 \left(4 \nu ^2-1\right) p^2+24\right) \\
& \quad -2 (n-7) (n-6) (n-5) p^4 \left(48 \nu +(2 \nu -1) p^2 \left(8 (\nu +1) (2 \nu +1)+p^2\right)+24\right) \\
& \quad +4 (n-6) (n-5) p^2 \left((3 \nu +2) p^6+\left(-4 \left(4 \nu ^2+6 \nu +1\right) \nu ^2+6 \nu +5\right) p^4+\left(3-12 \nu ^2\right) p^2+8\right) \\
& \quad +2 (n-5) p^2 \left(32 \nu +\left(4 (\nu +2) \nu ^2+\nu -1\right) p^6+4 (3 \nu +1) p^4+12 (\nu +2) (2 \nu -1) (2 \nu +1) p^2+16\right) \\\
& \quad -p^8-2 \left(4 \nu ^2+2 \nu -3\right) p^6+\left(12 \nu  \left(4 \nu ^2+2 \nu -1\right)-11\right) p^4+4 \left(1-4 \nu ^2\right) p^2-4
\end{aligned}
\right)
}{\left(4 \nu ^2-1\right) n (n+1) (2 \nu +n) (2 \nu +n+1)},
\end{equation*}
\begin{equation*}
\beta_6(n) = - \frac{p^4 \left(
\begin{aligned}
& +4 (n-9) (n-8) (n-7) (n-6) p^4 +2 (n-8) (n-7) (n-6) p^2 \left((4 \nu  (\nu +2)+9) p^2+12\right) \\
& \quad +6 (n-7) (n-6) p^2 \left(12 \nu +(2 \nu +1) \left(2 \nu ^2+\nu +2\right) p^2-2\right) \\
& \quad +2 (n-6) \left(\left((1-2 \nu )^2 \nu  (2 \nu +1)-1\right) p^4+(12 (\nu -4) \nu -31) p^2-12\right) \\
& \quad +(1-2 \nu ) p^4-2 \nu  (2 \nu +1) (6 \nu +1) p^2+8-24 \nu
\end{aligned}
\right)
}{\left(4 \nu ^2-1\right) n (n+1) (2 \nu +n) (2 \nu +n+1)},
\end{equation*}
\begin{equation*}
\beta_7(n) = - \frac{p^2 \left(
\begin{aligned}
& (n-10) (n-9) (n-8) (n-7) p^4 \left(\left(4 \nu ^2-1\right) p^2+16\right) \\
& \quad 4 (n-9) (n-8) (n-7) (2 \nu +1) p^4 \left(\left(2 \nu ^2+\nu -1\right) p^2+8\right) \\
& \quad +2 (n-8) (n-7) p^2 \left(\left(\nu  \left(2 \nu  \left(4 \nu ^2+6 \nu +1\right)-3\right)-4\right) p^4+4 \left(4 \nu ^2-1\right) p^2-24\right) \\
& \quad -2 (n-7) p^2 \left(48 \nu +(6 \nu +1) p^4+8 (\nu +2) (2 \nu -1) (2 \nu +1) p^2+24\right) \\
& \quad +(4 \nu  (\nu +1)+3) p^6+2 \left(7-4 \nu  \left(4 \nu ^2+2 \nu -1\right)\right) p^4+6 \left(4 \nu ^2-1\right) p^2+16
\end{aligned}
\right)
}{\left(4 \nu ^2-1\right) n (n+1) (2 \nu +n) (2 \nu +n+1)},
\end{equation*}
\begin{equation*}
\beta_8(n) = \frac{2 p^6\left(
\begin{aligned}
& (n-1) \left((4 (\nu -51) \nu +835) p^2-12\right) \\
& \quad -12 \nu +6 (2 \nu -17) (n-1)^2 p^2+4 (n-1)^3 p^2-(\nu  (4 \nu  (\nu +7)-835)+2222) p^2+98
\end{aligned}
\right)
}{\left(4 \nu ^2-1\right) n (n+1) (2 \nu +n) (2 \nu +n+1)},
\end{equation*}
\begin{equation*}
\beta_9(n) = \frac{2 p^4 \left(
\begin{aligned}
& 2 n^4+8 (\nu -10) n^3+\left(4 \nu ^2-240 \nu +1197\right) n^2 \\
& \quad -\left(8 \nu ^3+92 \nu ^2-2394 \nu +7937\right) n \\
& \quad -2 p^2 \left(-4 \nu ^2-144 \nu +8 n^2+16 (\nu -9) n+649\right)+p^4 \left(68 \nu ^3+502 \nu ^2-7937 \nu +19677\right)+12
\end{aligned}
\right)
}{\left(4 \nu ^2-1\right) n (n+1) (2 \nu +n) (2 \nu +n+1)},
\end{equation*}
\begin{equation*}
\beta_{10}(n) = \frac{8 p^8 (\nu +n-12)}{\left(4 \nu ^2-1\right) n (n+1) (2 \nu +n) (2 \nu +n+1)},
\end{equation*}
\begin{equation*}
\beta_{11}(n) = \frac{p^6 \left(p^2 \left(-4 \nu ^2-176 \nu +8 n^2+16 (\nu -11) n+969\right)-16\right)}{\left(4 \nu ^2-1\right) n (n+1) (2 \nu +n) (2 \nu +n+1)},
\end{equation*}
\begin{equation*}
\beta_{12}(n) = 0,
\end{equation*}
\begin{equation*}
\beta_{13}(n) = \frac{4 p^8}{\left(4 \nu ^2-1\right) n (n+1) (2 \nu +n) (2 \nu +n+1)}.
\end{equation*}
\end{thm}

\begin{thm}
Let $\nu,p \in \mathbb{C}$ and $\nu \notin \{-\frac{1}{2},-\frac{3}{2},\cdots\}$. Suppose that
\begin{equation*}
\arccos(pz) J_\nu(z) = \frac{1}{2^\nu \Gamma (\nu +1)} \sum_{n=0}^\infty u_n z^{n+\nu}, \quad z \in \mathbb{C}.
\end{equation*}
Then
\begin{equation*}
u_0 = \frac{\pi }{2}, u_1 = -p, u_2= -\frac{\pi }{8 (\nu +1)}, u_3 = \frac{p}{4 (\nu +1)}-\frac{p^3}{6}, u_4 = \frac{\pi }{64 (\nu +1) (\nu +2)},
\end{equation*}
\begin{equation*}
u_5 = -\frac{3 p^5}{40}+\frac{p^3}{24 (\nu +1)}-\frac{p}{32 (\nu +1) (\nu +2)}, u_6 = -\frac{\pi }{768 (\nu +1) (\nu +2) (\nu +3)},
\end{equation*}
\begin{equation*}
u_7 = -\frac{5 p^7}{112}+\frac{3 p^5}{160 (\nu +1)}-\frac{p^3}{192 (\nu +1) (\nu +2)}+\frac{p}{384 (\nu +1) (\nu +2) (\nu +3)},
\end{equation*}
\begin{equation*}
u_8 = \frac{\pi }{12288 (\nu +1) (\nu +2) (\nu +3) (\nu +4)},
\end{equation*}
\begin{equation*}
u_9 = -\frac{35 p^9}{1152}+\frac{5 p^7}{448 (\nu +1)}-\frac{3 p^5}{1280 (\nu +1) (\nu +2)}+\frac{p^3}{2304 (\nu +1) (\nu +2) (\nu +3)}-\frac{p}{6144 (\nu +1) (\nu +2) (\nu +3) (\nu +4)},
\end{equation*}
\begin{equation*}
u_{10} = -\frac{\pi }{245760 (\nu +1) (\nu +2) (\nu +3) (\nu +4) (\nu +5)},
\end{equation*}
\begin{equation*}
\begin{aligned}
u_{11} &= -\frac{63 p^{11}}{2816}+\frac{35 p^9}{4608 (\nu +1)}-\frac{5 p^7}{3584 (\nu +1) (\nu +2)}+\frac{p^5}{5120 (\nu +1) (\nu +2) (\nu +3)} \\
& \quad -\frac{p^3}{36864 (\nu +1) (\nu +2) (\nu +3) (\nu +4)}+\frac{p}{122880 (\nu +1) (\nu +2) (\nu +3) (\nu +4) (\nu +5)},
\end{aligned}
\end{equation*}
\begin{equation*}
u_{12} = \frac{\pi }{5898240 (\nu +1) (\nu +2) (\nu +3) (\nu +4) (\nu +5) (\nu +6)},
\end{equation*}
\begin{equation*}
\begin{aligned}
u_{13} &= -\frac{231 p^{13}}{13312}+\frac{63 p^{11}}{11264 (\nu +1)}-\frac{35 p^9}{36864 (\nu +1) (\nu +2)}+\frac{5 p^7}{43008 (\nu +1) (\nu +2) (\nu +3)}\\
& \quad -\frac{p^5}{81920 (\nu +1) (\nu +2) (\nu +3) (\nu +4)}+\frac{p^3}{737280 (\nu +1) (\nu +2) (\nu +3) (\nu +4) (\nu +5)}\\
& \quad -\frac{p}{2949120 (\nu +1) (\nu +2) (\nu +3) (\nu +4) (\nu +5) (\nu +6)},
\end{aligned}
\end{equation*}
and
\begin{equation*}
u_{n+1} = \sum_{i=0}^{13} \beta_i(n) u_{n-i}, n\geq13,
\end{equation*}
where $\beta_i(n)$ ($n=0,1,\cdots,13$) is defined in Theorem \ref{thm-arcsin-J}.
\end{thm}

\subsection{Recurrence Relations for the Maclaurin Coefficients of Products of Elementary Functions and the modified Bessel Function of the first kind}

We first consider the product of the exponential function and the modified Bessel function of the first kind.

\begin{thm} \label{thm-exp-I}
Let $\nu,p \in \mathbb{C}$ and $\nu \notin \{-\frac{1}{2},-\frac{3}{2},\cdots\}$. Suppose that
\begin{equation*}
e^{pz} I_\nu(z) = \frac{1}{2^\nu \Gamma (\nu +1)} \sum_{n=0}^\infty u_n z^{n+\nu}, \quad z \in \mathbb{C}.
\end{equation*}
Then $u_0 = 1$, $u_1 = p$ and
\begin{equation*}
u_{n+1} = \frac{p (2 \nu +2 n+1)}{(n+1) (2 \nu +n+1)} u_n + \frac{1-p^2}{(n+1) (2 \nu +n+1)} u_{n-1}.
\end{equation*}
\end{thm}

\begin{remark}
For $\nu,p\in\mathbb{R}$ and $z>0$, an equivalent recurrence was given by Yang \cite[Theorem~3.2]{Yang2025Report} in an internal seminar report.
\end{remark}

Noting that
\begin{equation*}
\sinh(pz) I_\nu(z) = \frac{1}{2} \Big( e^{pz} I_\nu(z) - e^{-pz} I_\nu(z) \Big), \end{equation*}
and
\begin{equation*}
\cosh(pz) I_\nu(z) = \frac{1}{2} \Big( e^{pz} I_\nu(z) + e^{-pz} I_\nu(z) \Big),
\end{equation*}
which lead to the following two theorems from Theorem \ref{thm-exp-I}.

\begin{thm} \label{thm-sinh-I-2}
Let $\nu,p \in \mathbb{C}$ and $\nu \notin \{-\frac{1}{2},-\frac{3}{2},\cdots\}$. Suppose that
\begin{equation*}
\sinh(pz) I_\nu(z) = \frac{1}{2^{\nu+1} \Gamma (\nu +1)} \sum_{n=0}^\infty (u_n - v_n) z^{n+\nu}, \quad z \in \mathbb{C}.
\end{equation*}
Then $u_0 = 1$, $u_1 = p$, $v_0 = 1$, $v_1 = -p$,
\begin{equation*}
u_{n+1} = \frac{p (2 \nu +2 n+1)}{(n+1) (2 \nu +n+1)} u_n + \frac{1-p^2}{(n+1) (2 \nu +n+1)} u_{n-1},
\end{equation*}
and
\begin{equation*}
v_{n+1} = - \frac{p (2 \nu +2 n+1)}{(n+1) (2 \nu +n+1)} v_n + \frac{1-p^2}{(n+1) (2 \nu +n+1)} v_{n-1}.
\end{equation*}
\end{thm}

\begin{thm} \label{thm-cosh-I-2}
Let $\nu,p \in \mathbb{C}$ and $\nu \notin \{-\frac{1}{2},-\frac{3}{2},\cdots\}$. Suppose that
\begin{equation*}
\cosh(pz) I_\nu(z) = \frac{1}{2^{\nu+1} \Gamma (\nu +1)} \sum_{n=0}^\infty (u_n + v_n) z^{n+\nu}, \quad z \in \mathbb{C}.
\end{equation*}
Then $u_0 = 1$, $u_1 = p$, $v_0 = 1$, $v_1 = -p$,
\begin{equation*}
u_{n+1} = \frac{p (2 \nu +2 n+1)}{(n+1) (2 \nu +n+1)} u_n + \frac{1-p^2}{(n+1) (2 \nu +n+1)} u_{n-1},
\end{equation*}
and
\begin{equation*}
v_{n+1} = - \frac{p (2 \nu +2 n+1)}{(n+1) (2 \nu +n+1)} v_n + \frac{1-p^2}{(n+1) (2 \nu +n+1)} v_{n-1}.
\end{equation*}
\end{thm}

Noting that
\begin{equation*}
\sin(pz) I_\nu(z) = \frac{1}{2i} \Big( e^{ipz} I_\nu(z) - e^{-ipz} I_\nu(z) \Big), \end{equation*}
and
\begin{equation*}
\cos(pz) I_\nu(z) = \frac{1}{2} \Big( e^{ipz} I_\nu(z) + e^{-ipz} I_\nu(z) \Big),
\end{equation*}
which lead to the following two theorems from Theorem \ref{thm-exp-I}.

\begin{thm} \label{thm-sin-I-2}
Let $\nu,p \in \mathbb{C}$ and $\nu \notin \{-\frac{1}{2},-\frac{3}{2},\cdots\}$. Suppose that
\begin{equation*}
\sin(pz) I_\nu(z) = \frac{1}{2^{\nu+1} i \Gamma (\nu +1)} \sum_{n=0}^\infty (u_n - v_n) z^{n+\nu}, \quad z \in \mathbb{C}.
\end{equation*}
Then $u_0 = 1$, $u_1 = ip$, $v_0 = 1$, $v_1 = -ip$,
\begin{equation*}
u_{n+1} = \frac{ip (2 \nu +2 n+1)}{(n+1) (2 \nu +n+1)} u_n + \frac{p^2+1}{(n+1) (2 \nu +n+1)} u_{n-1},
\end{equation*}
and
\begin{equation*}
v_{n+1} = - \frac{i p (2 \nu +2 n+1)}{(n+1) (2 \nu +n+1)} v_n + \frac{p^2+1}{(n+1) (2 \nu +n+1)} v_{n-1}.
\end{equation*}
\end{thm}

\begin{thm} \label{thm-cos-I-2}
Let $\nu,p \in \mathbb{C}$ and $\nu \notin \{-\frac{1}{2},-\frac{3}{2},\cdots\}$. Suppose that
\begin{equation*}
\cos(pz) I_\nu(z) = \frac{1}{2^{\nu+1} \Gamma (\nu +1)} \sum_{n=0}^\infty (u_n + v_n) z^{n+\nu}, \quad z \in \mathbb{C}.
\end{equation*}
Then $u_0 = 1$, $u_1 = ip$, $v_0 = 1$, $v_1 = -ip$,
\begin{equation*}
u_{n+1} = \frac{ip (2 \nu +2 n+1)}{(n+1) (2 \nu +n+1)} u_n + \frac{p^2+1}{(n+1) (2 \nu +n+1)} u_{n-1},
\end{equation*}
and
\begin{equation*}
v_{n+1} = - \frac{i p (2 \nu +2 n+1)}{(n+1) (2 \nu +n+1)} v_n + \frac{p^2+1}{(n+1) (2 \nu +n+1)} v_{n-1}.
\end{equation*}
\end{thm}

Next, we consider the product of $(1-\theta z)^p$ and the modified Bessel function of the first kind.

\begin{thm}
Let $\nu,p,\theta \in \mathbb{C}$ and $\nu \notin \{-\frac{1}{2},-\frac{3}{2},\cdots\}$. Suppose that
\begin{equation*}
(1-\theta z)^p I_\nu(z) = \frac{1}{2^\nu \Gamma (\nu +1)} \sum_{n=0}^\infty u_n z^{n+\nu}, \quad |z| < \frac{1}{\theta}.
\end{equation*}
Then $u_0 = 1$, $u_1 = -\theta p$, $u_2 = \frac{1}{2} \theta ^2 (p-1) p + \frac{1}{4 (\nu +1)}$,
$u_3 = -\frac{1}{6} \theta ^3 (p-2) (p-1) p-\frac{\theta  p}{4 (\nu +1)}$ and
\begin{equation*}
u_{n+1} = \beta_0(n) u_n + \beta_1(n) u_{n-1} + \beta_2(n) u_{n-2} + \beta_3(n) u_{n-3},
\end{equation*}
where
\begin{equation*}
\begin{aligned}
\beta_0(n) &= \frac{\theta  \left(2 n^2+4 \nu  n-2 n p-2 \nu  p-p\right)}{(n+1) (2 \nu +n+1)},\\
\beta_1(n) &= \frac{1-\theta ^2 (n-p-1) (2 \nu +n-p-1)}{(n+1) (2 \nu +n+1)},\\
\beta_2(n) &= -\frac{2 \theta }{(n+1) (2 \nu +n+1)}, \\
\beta_3(n) &= \frac{\theta ^2}{(n+1) (2 \nu +n+1)}.
\end{aligned}
\end{equation*}
\end{thm}

Next, we consider the product of $\exp(-p \arctan z)$ and the modified Bessel function of the first kind.

\begin{thm}
Let $\nu,p \in \mathbb{C}$ and $\nu \notin \{-\frac{1}{2},-\frac{3}{2},\cdots\}$. Suppose that
\begin{equation*}
e^{-p \arctan z} I_\nu(z) = \frac{1}{2^\nu \Gamma (\nu +1)} \sum_{n=0}^\infty u_n z^{n+\nu}, \quad z \in \mathbb{C}.
\end{equation*}
Then $u_0 = 1$, $u_1 = - p$, $u_2 = \frac{p^2}{2}+\frac{1}{4 (\nu +1)}$,
$u_3 = \frac{1}{3} (p-\frac{p^3}{2})-\frac{p}{4 (\nu +1)}$,
$u_4 = \frac{1}{96} (\frac{3}{\nu ^2+3 \nu +2}+4 p^2 (\frac{3}{\nu +1}+p^2-8))$, and
\begin{equation*}
u_{n+1} = \beta_0(n) u_n + \beta_1(n) u_{n-1} + \beta_2(n) u_{n-2} + \beta_3(n) u_{n-3} + \beta_4(n) u_{n-4},
\end{equation*}
where
\begin{equation*}
\begin{aligned}
\beta_0(n) &= -\frac{p (2 \nu +2 n+1)}{(n+1) (2 \nu +n+1)},\\
\beta_1(n) &= -\frac{-4 \nu +2 n (2 \nu +n-2)+p^2+1}{(n+1) (2 \nu +n+1)},\\
\beta_2(n) &= \frac{p (-2 \nu -2 n+5)}{(n+1) (2 \nu +n+1)}, \\
\beta_3(n) &= \frac{6 \nu -n (2 \nu +n-6)-7}{(n+1) (2 \nu +n+1)}, \\
\beta_4(n) &= \frac{1}{(n+1) (2 \nu +n+1)}.
\end{aligned}
\end{equation*}
\end{thm}

Next, the following two theorems consider the products of trigonometric functions and the modified Bessel function of the first kind, and characterize their coefficients by a single recurrence relation.

\begin{thm} \label{thm-sin-I}
Let $\nu,p \in \mathbb{C}$ and $\nu \notin \{-\frac{1}{2},-\frac{3}{2},\cdots\}$. Suppose that
\begin{equation*}
\sin(pz) I_\nu(z) = \frac{1}{2^\nu \Gamma (\nu +1)} \sum_{n=0}^\infty u_n z^{n+\nu}, \quad z \in \mathbb{C}.
\end{equation*}
Then $u_{2k}=0$ for all $k\geq 0$, $u_1 = p$, $u_3 = -\frac{p^3}{6}+\frac{p}{4 (\nu +1)}$, $u_5 = \frac{p^5}{120}-\frac{p^3}{24 (\nu +1)}+\frac{p}{32 (\nu +1) (\nu +2)}$ and
\begin{equation*}
u_{2n+2} = \beta_0(n) u_{2n} + \beta_1(n) u_{2n-2} + \beta_2(n) u_{2n-4},
\end{equation*}
where
\begin{equation*}
\begin{aligned}
\beta_0(n) &= -\frac{\left(
\begin{aligned}
& 4 \left(p^2+1\right)n^4 +16 (\nu -1) \left(p^2+1\right)n^3 +2 \left(4 (\nu -6) \nu +(12 (\nu -2) \nu +7) p^2+9\right)n^2 \\
& \quad +2 \left(-2 \nu  (2 \nu  (2 \nu +5)-9)+\left(2 \nu  \left(4 \nu ^2-2 \nu +7\right)-3\right) p^2+1\right) n \\
& \quad +(2 \nu -1) (2 \nu +1)^2 \left(2 (\nu +1) p^2-1\right)
\end{aligned}
\right)}{\left(4 \nu ^2-1\right) (n+1) (n+2) (2 \nu +n+1) (2 \nu +n+2)},\\
\beta_1(n) &= -\frac{8 n^2 \left(p^4-1\right)+16 (\nu -2) n \left(p^4-1\right)+\left(p^2+1\right) \left(4 \nu  (\nu +8)+(4 \nu  (5 \nu -8)+27) p^2-33\right)}{\left(4 \nu ^2-1\right) (n+1) (n+2) (2 \nu +n+1) (2 \nu +n+2)},\\
\beta_2(n) &= -\frac{4 \left(p^2+1\right)^3}{\left(4 \nu ^2-1\right) (n+1) (n+2) (2 \nu +n+1) (2 \nu +n+2)}.
\end{aligned}
\end{equation*}
\end{thm}

\begin{thm}
Let $\nu,p \in \mathbb{C}$ and $\nu \notin \{-\frac{1}{2},-\frac{3}{2},\cdots\}$. Suppose that
\begin{equation*}
\cos(pz) I_\nu(z) = \frac{1}{2^\nu \Gamma (\nu +1)} \sum_{n=0}^\infty u_n z^{n+\nu}, \quad z \in \mathbb{C}.
\end{equation*}
Then $u_{2k-1}=0$ for all $k\geq 0$, $u_0 = 1$, $u_2 =\frac{1}{4 (\nu +1)}-\frac{p^2}{2}$, $u_4 = \frac{1}{32 (\nu +1) (\nu +2)}+\frac{p^4}{24}-\frac{p^2}{8 (\nu +1)}$ and
\begin{equation*}
u_{2n+2} = \beta_0(n) u_{2n} + \beta_1(n) u_{2n-2} + \beta_2(n) u_{2n-4},
\end{equation*}
where $\beta_i(n)$ ($n=0,1,2$) is defined in Theorem \ref{thm-sin-I}.
\end{thm}

Next, we consider the product of hyperbolic functions and the modified Bessel function of the first kind.

\begin{thm} \label{thm-sinh-I}
Let $\nu,p \in \mathbb{C}$ and $\nu \notin \{-\frac{1}{2},-\frac{3}{2},\cdots\}$. Suppose that
\begin{equation*}
\sinh(pz) I_\nu(z) = \frac{1}{2^\nu \Gamma (\nu +1)} \sum_{n=0}^\infty u_n z^{n+\nu}, \quad z \in \mathbb{C}.
\end{equation*}
Then $u_{2k}=0$ for all $k\geq 0$, $u_1 = p$, $u_3 = \frac{p^3}{6}+\frac{p}{4 (\nu +1)}$, $u_5 = \frac{p^5}{120}+\frac{p^3}{24 (\nu +1)}+\frac{p}{32 (\nu +1) (\nu +2)}$ and
\begin{equation*}
u_{2n+2} = \beta_0(n) u_{2n} + \beta_1(n) u_{2n-2} + \beta_2(n) u_{2n-4},
\end{equation*}
where
\begin{equation*}
\begin{aligned}
\beta_0(n) &= \frac{
\left(
\begin{aligned}
& 4 \left(p^2-1\right)n^4 +16 (\nu -1) \left(p^2-1\right)n^3+2 \left(-4 (\nu -6) \nu +(12 (\nu -2) \nu +7) p^2-9\right)n^2 \\
& \quad +2 \left(2 \nu  (2 \nu  (2 \nu +5)-9)+\left(2 \nu  \left(4 \nu ^2-2 \nu +7\right)-3\right) p^2-1\right)n \\
& \quad +(2 \nu -1) (2 \nu +1)^2 \left(2 (\nu +1) p^2+1\right)
\end{aligned}
\right)}{\left(4 \nu ^2-1\right) (n+1) (n+2) (2 \nu +n+1) (2 \nu +n+2)},\\
\beta_1(n) &= -\frac{8 n^2 \left(p^4-1\right)+16 (\nu -2) n \left(p^4-1\right)+\left(p^2-1\right) \left(-4 \nu  (\nu +8)+(4 \nu  (5 \nu -8)+27) p^2+33\right)}{\left(4 \nu ^2-1\right) (n+1) (n+2) (2 \nu +n+1) (2 \nu +n+2)},\\
\beta_2(n) &= \frac{4 \left(p^2-1\right)^3}{\left(4 \nu ^2-1\right) (n+1) (n+2) (2 \nu +n+1) (2 \nu +n+2)}.
\end{aligned}
\end{equation*}
\end{thm}

\begin{thm}
Let $\nu,p \in \mathbb{C}$ and $\nu \notin \{-\frac{1}{2},-\frac{3}{2},\cdots\}$. Suppose that
\begin{equation*}
\cosh(pz) I_\nu(z) = \frac{1}{2^\nu \Gamma (\nu +1)} \sum_{n=0}^\infty u_n z^{n+\nu}, \quad z \in \mathbb{C}.
\end{equation*}
Then $u_{2k-1}=0$ for all $k\geq 0$, $u_0 = 1$, $u_2 = \frac{p^2}{2}+\frac{1}{4 (\nu +1)}$, $u_4 = \frac{1}{96} (-\frac{3}{\nu +2}+4 p^4+\frac{12 p^2+3}{\nu +1})$ and
\begin{equation*}
u_{2n+2} = \beta_0(n) u_{2n} + \beta_1(n) u_{2n-2} + \beta_2(n) u_{2n-4},
\end{equation*}
where $\beta_i(n)$ ($n=0,1,2$) is defined in Theorem \ref{thm-sinh-I}.
\end{thm}

At the end, we consider the product of an inverse trigonometric function and a modified Bessel function of the first kind.

\begin{thm} \label{thm-arcsin-I}
Let $\nu,p \in \mathbb{C}$ and $\nu \notin \{-\frac{1}{2},-\frac{3}{2},\cdots\}$. Suppose that
\begin{equation*}
\arcsin(pz) I_\nu(z) = \frac{1}{2^\nu \Gamma (\nu +1)} \sum_{n=0}^\infty u_n z^{n+\nu}, \quad z \in \mathbb{C}.
\end{equation*}
Then
$u_0 = u_2 = u_4 = u_6 = u_8 = u_{10} = u_{12} = 0$,
\begin{equation*}
u_1 = p, u_3 = \frac{p^3}{6}+\frac{p}{4 (\nu +1)}, u_5 = \frac{3 p^5}{40}+\frac{p^3}{24 (\nu +1)}+\frac{p}{32 (\nu +1) (\nu +2)},
\end{equation*}
\begin{equation*}
u_7 = \frac{5 p^7}{112}+\frac{3 p^5}{160 (\nu +1)}+\frac{p^3}{192 (\nu +1) (\nu +2)}+\frac{p}{384 (\nu +1) (\nu +2) (\nu +3)},
\end{equation*}
\begin{equation*}
u_9 = \frac{35 p^9}{1152}+\frac{5 p^7}{448 (\nu +1)}+\frac{3 p^5}{1280 (\nu +1) (\nu +2)}+\frac{p^3}{2304 (\nu +1) (\nu +2) (\nu +3)}+\frac{p}{6144 (\nu +1) (\nu +2) (\nu +3) (\nu +4)},
\end{equation*}
\begin{equation*}
\begin{aligned}
u_{11} & = \frac{63 p^{11}}{2816}+\frac{35 p^9}{4608 (\nu +1)}+\frac{5 p^7}{3584 (\nu +1) (\nu +2)}+\frac{p^5}{5120 (\nu +1) (\nu +2) (\nu +3)} \\
& \quad +\frac{p^3}{36864 (\nu +1) (\nu +2) (\nu +3) (\nu +4)}+\frac{p}{122880 (\nu +1) (\nu +2) (\nu +3) (\nu +4) (\nu +5)},
\end{aligned}
\end{equation*}
\begin{equation*}
\begin{aligned}
u_{13} & = \frac{231 p^{13}}{13312}+\frac{63 p^{11}}{11264 (\nu +1)}+\frac{35 p^9}{36864 (\nu +1) (\nu +2)}+\frac{5 p^7}{43008 (\nu +1) (\nu +2) (\nu +3)}\\
& \quad +\frac{p^5}{81920 (\nu +1) (\nu +2) (\nu +3) (\nu +4)}+\frac{p^3}{737280 (\nu +1) (\nu +2) (\nu +3) (\nu +4) (\nu +5)}\\
& \quad +\frac{p}{2949120 (\nu +1) (\nu +2) (\nu +3) (\nu +4) (\nu +5) (\nu +6)},
\end{aligned}
\end{equation*}
and
\begin{equation*}
u_{n+1} = \sum_{i=0}^{13} \beta_i(n) u_{n-i}, n\geq 13,
\end{equation*}
where
\begin{equation*}
\beta_0(n) = \frac{2 p^2 (\nu +n)}{(n+1) (2 \nu +n+1)},
\end{equation*}
\begin{equation*}
\beta_1(n) = \frac{\left(
\begin{aligned}
& (n-4) (n-3) (n-2) (n-1) \left(p^4+4 \left(4 \nu ^2-1\right) p^2-4\right) \\
& \quad +2 (2 \nu +1) (n-3) (n-2) (n-1) \left(p^4+8 \left(2 \nu ^2+\nu -1\right) p^2-4\right) \\
& \quad +2 \left(4 \nu ^2-1\right) (n-2) (n-1) \left(4 (\nu +1) (2 \nu +1) p^2-1\right) \\
& \quad (2 \nu -1) (2 \nu +1)^2-2 (2 \nu -1) (2 \nu +1) (n-1) \left((\nu +1) p^4-2 (\nu +2)\right)
\end{aligned}
\right)
}{(4 \nu ^2-1) n (n+1) (2 \nu +n) (2 \nu +n+1)},
\end{equation*}
\begin{equation*}
\beta_2(n) = - \frac{p^2 \left(
\begin{aligned}
& 4 (n-5) (n-4) (n-3) (n-2) p^2 \\
& \quad + 2 (n-4) (n-3) (n-2) \left(p^4+(4 \nu  (3 \nu +2)-1) p^2-4\right) \\
& \quad + 3 (2 \nu +1) (n-3) (n-2) \left(p^2 \left(6 \nu  (2 \nu +1)+p^2-8\right)-4\right) \\
& \quad + \left(4 \nu ^2-1\right) (n-2) \left(p^4+2 (\nu  (6 \nu +5)-2) p^2-2\right) \\
& \quad + 2 (\nu +2) (2 \nu -1) (2 \nu +1)
\end{aligned}
\right)
}{(4 \nu ^2-1) n (n+1) (2 \nu +n) (2 \nu +n+1)},
\end{equation*}
\begin{equation*}
\beta_3(n) = \frac{\left(
\begin{aligned}
& -2 (n-6) (n-5) (n-4) (n-3) p^2 \left(p^4+3 \left(4 \nu ^2-1\right) p^2-8\right) \\
& \quad -8 (n-5) (n-4) (n-3) p^2 \left(-8 \nu +p^2 \left(\nu  \left(12 \nu  (\nu +1)+p^2-3\right)-3\right)-4\right) \\
& \quad +(n-4) (n-3) \left(p^8+12 \nu  p^6-6 \left(2 \nu  \left(2 \nu  \left(4 \nu ^2+6 \nu +1\right)-3\right)-1\right) p^4+8 \left(4 \nu ^2-1\right) p^2+8\right) \\
& \quad +(2 \nu +1) (n-3) \left(p^8+4 (2 \nu  (\nu +1)-3) p^6-6 p^4-16 (\nu +2) (2 \nu -1) p^2+8\right) \\
& \quad -(2 \nu -1) (2 \nu +1) \left(-p^4+(8 \nu +4) p^2+1\right)
\end{aligned}
\right)
}{(4 \nu ^2-1) n (n+1) (2 \nu +n) (2 \nu +n+1)},
\end{equation*}
\begin{equation*}
\beta_4(n) = \frac{p^2 \left(
\begin{aligned}
& 8 n^4 p^4 +2 n^3 p^2 \left(p^4+\left(12 \nu ^2+16 \nu -79\right) p^2-12\right) \\
& \quad +n^2 p^2 \left(-72 \nu +(6 \nu -29) p^4+2 \left(36 \nu ^3-144 \nu ^2-237 \nu +578\right) p^2+348\right) \\
& \quad +2 n \left(\left(2 \nu ^2-29 \nu +68\right) p^6+\left(24 \nu ^4-320 \nu ^3+546 \nu ^2+1156 \nu -1855\right) p^4+\left(-12 \nu ^2+348 \nu -821\right) p^2-4\right) \\
& \quad -8 \nu -8 \left(2 \nu ^2-17 \nu +26\right) p^6+\left(-192 \nu ^4+1408 \nu ^3-1296 \nu ^2-3710 \nu +4409\right) p^4 \\
& \quad +2 \left(12 \nu ^3+64 \nu ^2-821 \nu +1263\right) p^2+28
\end{aligned}
\right)
}{(4 \nu ^2-1) n (n+1) (2 \nu +n) (2 \nu +n+1)},
\end{equation*}
\begin{equation*}
\beta_5(n) = \frac{\left(
\begin{aligned}
& (n-8) (n-7) (n-6) (n-5) p^4 \left(p^4+4 \left(4 \nu ^2-1\right) p^2-24\right) \\
& \quad +2 (n-7) (n-6) (n-5) \left((2 \nu -1) p^8+8 (\nu +1) (2 \nu -1) (2 \nu +1) p^6-24 (2 \nu +1) p^4\right) \\
& \quad -4 (n-6) (n-5) p^2 \left((3 \nu +2) p^6-\left(2 \nu  \left(2 \nu  \left(4 \nu ^2+6 \nu +1\right)-3\right)+1\right) p^4+3 \left(4 \nu ^2-1\right) p^2+8\right) \\
& \quad -2 (n-5) p^2 \left(32 \nu +\left(4 (\nu +2) \nu ^2+\nu -1\right) p^6-4 (3 \nu +1) p^4+12 \left(-4 \nu ^3-8 \nu ^2+\nu +2\right) p^2+16\right) \\
& \quad -p^8-2 \left(4 \nu ^2+2 \nu -3\right) p^6+\left(12 \nu  \left(4 \nu ^2+2 \nu -1\right)-1\right) p^4+4 \left(4 \nu ^2-1\right) p^2-4
\end{aligned}
\right)
}{(4 \nu ^2-1) n (n+1) (2 \nu +n) (2 \nu +n+1)},
\end{equation*}
\begin{equation*}
\beta_6(n) = \frac{p^4 \left(
\begin{aligned}
& -4 (n-9) (n-8) (n-7) (n-6) p^4 -2 (n-8) (n-7) (n-6) p^2 \left((4 \nu  (\nu +2)+9) p^2-12\right) \\
& \quad -6 (n-7) (n-6) p^2 \left(-12 \nu +(2 \nu +1) \left(2 \nu ^2+\nu +2\right) p^2+2\right) \\
& \quad -2 (n-6) \left(\left(\nu  (2 \nu +1) (1-2 \nu )^2+1\right) p^4+(31-12 (\nu -4) \nu ) p^2-12\right) \\
& \quad 24 \nu +(1-2 \nu ) p^4-2 \nu  (2 \nu +1) (6 \nu +1) p^2-8
\end{aligned}
\right)
}{(4 \nu ^2-1) n (n+1) (2 \nu +n) (2 \nu +n+1)},
\end{equation*}
\begin{equation*}
\beta_7(n) = - \frac{\left(
\begin{aligned}
& (n-10) (n-9) (n-8) (n-7) p^6 \left(\left(4 \nu ^2-1\right) p^2-16\right) \\
& \quad + 4 (2 \nu +1) (n-9) (n-8) (n-7) p^6 \left(\left(2 \nu ^2+\nu -1\right) p^2-8\right) \\
& \quad +2 (n-8) (n-7) p^4 \left(\left(\nu  \left(2 \nu  \left(4 \nu ^2+6 \nu +1\right)-3\right)+2\right) p^4+\left(4-16 \nu ^2\right) p^2-24\right) \\
& \quad +2 (n-7) \left((6 \nu +1) p^8+8 (\nu +2) (2 \nu -1) (2 \nu +1) p^6-24 (2 \nu +1) p^4\right) \\
& \quad -p^2 \left((4 \nu  (\nu +1)+3) p^6-2 \left(4 \nu  \left(4 \nu ^2+2 \nu -1\right)+3\right) p^4+\left(6-24 \nu ^2\right) p^2+16\right)
\end{aligned}
\right)
}{(4 \nu ^2-1) n (n+1) (2 \nu +n) (2 \nu +n+1)},
\end{equation*}
\begin{equation*}
\beta_8(n) = \frac{2 p^6\left(
\begin{aligned}
& -4 n^3 p^2 +6 (19-2 \nu ) n^2 p^2  -n \left((4 (\nu -57) \nu +1051) p^2+12\right) \\
& \quad -12 \nu +(\nu  (4 \nu  (\nu +8)-1051)+3163) p^2+110
\end{aligned}
\right)
}{(4 \nu ^2-1) n (n+1) (2 \nu +n) (2 \nu +n+1)},
\end{equation*}
\begin{equation*}
\beta_9(n) = \frac{2 p^4 \left(
\begin{aligned}
& -2 n^4 p^4 - 8 (\nu -10) n^3 p^4 -n^2 p^2 \left(\left(4 \nu ^2-240 \nu +1197\right) p^2+16\right) \\
& \quad +n \left(\left(8 \nu ^3+92 \nu ^2-2394 \nu +7937\right) p^4-32 (\nu -9) p^2\right) \\
& \quad -\left(68 \nu ^3+502 \nu ^2-7937 \nu +19672\right) p^4+2 \left(4 \nu ^2+144 \nu -649\right) p^2-12
\end{aligned}
\right)
}{(4 \nu ^2-1) n (n+1) (2 \nu +n) (2 \nu +n+1)},
\end{equation*}
\begin{equation*}
\beta_{10}(n) = \frac{8 p^8 (\nu +n-12)}{\left(4 \nu ^2-1\right) n (n+1) (2 \nu +n) (2 \nu +n+1)},
\end{equation*}
\begin{equation*}
\beta_{11}(n) = \frac{p^6 \left(p^2 \left(-4 \nu  (\nu +44)+8 n^2+16 (\nu -11) n+969\right)+16\right)}{\left(4 \nu ^2-1\right) n (n+1) (2 \nu +n) (2 \nu +n+1)},
\end{equation*}
\begin{equation*}
\beta_{12}(n) = 0,
\end{equation*}
\begin{equation*}
\beta_{13}(n) = -\frac{4 p^8}{\left(4 \nu ^2-1\right) n (n+1) (2 \nu +n) (2 \nu +n+1)}.
\end{equation*}
\end{thm}

\begin{thm}
Let $\nu,p \in \mathbb{C}$ and $\nu \notin \{-\frac{1}{2},-\frac{3}{2},\cdots\}$. Suppose that
\begin{equation*}
\arccos(pz) J_\nu(z) = \frac{1}{2^\nu \Gamma (\nu +1)} \sum_{n=0}^\infty u_n z^{n+\nu}, \quad z \in \mathbb{C}.
\end{equation*}
Then
\begin{equation*}
u_0 = \frac{\pi }{2}, u_1 = -p, u_2= \frac{\pi }{8 (\nu +1)}, u_3 = -\frac{p^3}{6}-\frac{p}{4 (\nu +1)}, u_4 = \frac{\pi }{64 (\nu +1) (\nu +2)},
\end{equation*}
\begin{equation*}
u_5 = -\frac{3 p^5}{40}-\frac{p^3}{24 (\nu +1)}-\frac{p}{32 (\nu +1) (\nu +2)}, u_6 = \frac{\pi }{768 (\nu +1) (\nu +2) (\nu +3)},
\end{equation*}
\begin{equation*}
u_7 = -\frac{5 p^7}{112}-\frac{3 p^5}{160 (\nu +1)}-\frac{p^3}{192 (\nu +1) (\nu +2)}-\frac{p}{384 (\nu +1) (\nu +2) (\nu +3)},
\end{equation*}
\begin{equation*}
u_8 = \frac{\pi }{12288 (\nu +1) (\nu +2) (\nu +3) (\nu +4)},
\end{equation*}
\begin{equation*}
u_9 = -\frac{35 p^9}{1152}-\frac{5 p^7}{448 (\nu +1)}-\frac{3 p^5}{1280 (\nu +1) (\nu +2)}-\frac{p^3}{2304 (\nu +1) (\nu +2) (\nu +3)}-\frac{p}{6144 (\nu +1) (\nu +2) (\nu +3) (\nu +4)},
\end{equation*}
\begin{equation*}
u_{10} = \frac{\pi }{245760 (\nu +1) (\nu +2) (\nu +3) (\nu +4) (\nu +5)},
\end{equation*}
\begin{equation*}
\begin{aligned}
u_{11} &= -\frac{63 p^{11}}{2816}-\frac{35 p^9}{4608 (\nu +1)}-\frac{5 p^7}{3584 (\nu +1) (\nu +2)}-\frac{p^5}{5120 (\nu +1) (\nu +2) (\nu +3)}\\
& \quad -\frac{p^3}{36864 (\nu +1) (\nu +2) (\nu +3) (\nu +4)}-\frac{p}{122880 (\nu +1) (\nu +2) (\nu +3) (\nu +4) (\nu +5)},
\end{aligned}
\end{equation*}
\begin{equation*}
u_{12} = \frac{\pi }{5898240 (\nu +1) (\nu +2) (\nu +3) (\nu +4) (\nu +5) (\nu +6)},
\end{equation*}
\begin{equation*}
\begin{aligned}
u_{13} &= -\frac{231 p^{13}}{13312}-\frac{63 p^{11}}{11264 (\nu +1)}-\frac{35 p^9}{36864 (\nu +1) (\nu +2)}-\frac{5 p^7}{43008 (\nu +1) (\nu +2) (\nu +3)} \\
& \quad -\frac{p^5}{81920 (\nu +1) (\nu +2) (\nu +3) (\nu +4)}-\frac{p^3}{737280 (\nu +1) (\nu +2) (\nu +3) (\nu +4) (\nu +5)}\\
& \quad -\frac{p}{2949120 (\nu +1) (\nu +2) (\nu +3) (\nu +4) (\nu +5) (\nu +6)},
\end{aligned}
\end{equation*}
and
\begin{equation*}
u_{n+1} = \sum_{i=0}^{13} \beta_i(n) u_{n-i}, n\geq13,
\end{equation*}
where $\beta_i(n)$ ($n=0,1,2,\cdots,13$) is defined in Theorem \ref{thm-arcsin-I}.
\end{thm}


\begin{thebibliography}{1}

\bibitem{Bowman-1958}
F.~Bowman.
\newblock {\em Introduction to {B}essel functions}.
\newblock Dover Publications, Inc., New York, 1958.

\bibitem{Chen-RACSAM-2021}
Y.-J. Chen and T.-H. Zhao.
\newblock On the monotonicity and convexity for generalized elliptic integral
  of the first kind.
\newblock {\em Rev. R. Acad. Cienc. Exactas F\'{\i}s. Nat. Ser. A Mat. RACSAM},
  116(2):1--21, Paper No. 77, 2022.

\bibitem{Korenev-2002}
B.~G. Korenev.
\newblock {\em Bessel functions and their applications}, volume~8 of {\em
  Analytical Methods and Special Functions}.
\newblock Taylor \& Francis Group, London, 2002.
\newblock Translated from the Russian by E. V. Pankratiev.

\bibitem{Luke-1962}
Y.~L. Luke.
\newblock {\em Integrals of {B}essel functions}.
\newblock McGraw-Hill Book Co., Inc., New York-Toronto-London, 1962.

\bibitem{Wu-IJPAM-2025}
J.~Wu and T.~Zhao.
\newblock On the power series related to zero-balanced hypergeometric function.
\newblock {\em Indian J. Pure Appl. Math.}, 2025.

\bibitem{Yang-PAMS-2025}
Z.-H. Yang.
\newblock Recurrence relations of {M}aclaurin series coefficients involving
  hypergeometric functions with applications.
\newblock {\em Proc. Amer. Math. Soc.}, 153(8):3513--3527, 2025.

\bibitem{Yang2025Report}
Z.-H. Yang.
\newblock {\em Recurrence method in certain special functions}.
\newblock Internal seminar report, Hangzhou Normal University,
  September 26, 2025.

\bibitem{Yang-JMAA-2018}
Z.-H. Yang, W.-M. Qian, Y.-M. Chu, and W.~Zhang.
\newblock On approximating the arithmetic-geometric mean and complete elliptic
  integral of the first kind.
\newblock {\em J. Math. Anal. Appl.}, 462(2):1714--1726, 2018.

\end{thebibliography}

\end{document}